\documentclass[11pt,a4paper]{article}
\usepackage{indentfirst}
\usepackage{amsmath}
\usepackage{amsthm}
\usepackage{amssymb}
\usepackage{epsfig}
\usepackage{color}
\IfFileExists{ifsym.sty}{\usepackage[misc]{ifsym}}{%
  }
\usepackage{hyperref}
\usepackage{enumitem}
\numberwithin{equation}{section}
\newtheorem{definition}{Definition}[section]
\newtheorem{proposition}[definition]{Proposition}
\newtheorem{corollary}[definition]{Corollary}
\newtheorem{theorem}[definition]{Theorem}
\newtheorem{lemma}[definition]{Lemma}

\newtheorem{claim}[definition]{Claim}
\newtheorem{remark}[definition]{Remark}

\def\be{\begin{equation}}
\def\ee{\end{equation}}
\def\bea{\begin{eqnarray}}
\def\eea{\end{eqnarray}}
\def\bes{\begin{eqnarray*}}
\def\ees{\end{eqnarray*}}

\def\<{\langle}
\def\>{\rangle}
\def\lb{\label}

\def\T{{\cal T}}
\def\H{{\cal H}}
\def\R{{\bf R}}
\def\C{{\bf C}}
\def\Z{{\bf Z}}
\def\N{{\bf N}}
\def\U{{\bf U}}
\def\Q{{\bf Q}}
\def\CP{{\bf CP}}

\def\Sg{{\Sigma}}

\newcommand{\Sp}{\operatorname{Sp}}
\newcommand{\per}{\operatorname{per}}
\newcommand{\cT}{\mathcal T}

\newcommand{\cN}{\mathcal N}

\title{Nonvanishing of mean indices and stability of symmetric closed characteristics
on symmetric star-shaped hypersurfaces in $\R^4$}
\author{Zhenxiong Li $^{1}$,\thanks{E-mail: zxli@gxu.edu.cn} \quad
Hui Liu $^{2}$,\thanks{Partially supported by NSFC (12371195 and 12022111) and the Fundamental Research Funds for the Central
Universities (2042023kf0207). E-mail: huiliu00031514@whu.edu.cn.}
\quad Zhoukai Xu $^{3}$,\thanks{E-mail: zhoukaixu@sdu.edu.cn}\\
$^{1}$ School of Mathematics, Guangxi University, Nanning 530004, \\
Guangxi, the People's Republic of China \\
$^{2}$ School of Mathematics and Statistics, Wuhan University, Wuhan 430072, \\
Hubei, the People's Republic of China \\
$^{3}$ School of Mathematics, Shandong University, Jinan 250100, \\
Shandong, the People's Republic of China \\ }
\date{}

\begin{document}
\maketitle
\begin{abstract}
H. Duan, H. Liu, Y. Long and W. Wang \cite{DLLW24} conjectured that if a compact star-shaped hypersurface in $\R^{2n}$ carries only finitely many geometrically distinct prime closed characteristics, then each of them has nonzero mean index.
They proved this conjecture for non-degenerate star-shaped hypersurfaces in $\R^6$.
Motivated by their work, we study the symmetric counterpart of this problem.
We prove that if a symmetric compact star-shaped hypersurface $\Sigma$  in $\R^4$ carries only finitely many geometrically distinct prime symmetric closed characteristics, and every odd iterate of them is non-degenerate in the symmetric sense, 
then each of them has nonzero symmetric mean index.
Based on this result, we further prove that if $\Sigma$ carries exactly two geometrically distinct prime symmetric closed characteristics and all their iterates are non-degenerate, then both of them are irrationally elliptic.
\end{abstract}



\section{Introduction and main results}\label{sec:1}

Let $\Sigma$ be a $C^3$ compact hypersurface in $\R^{2n}$ strictly star-shaped with respect to the origin, i.e., the tangent hyperplane at any $x\in\Sigma$ does not intersect the origin.
We denote the set of all such hypersurfaces by $\H_{st}(2n)$, and denote by $\H_{con}(2n)$ the subset of $\H_{st}(2n)$ which consists of all strictly convex hypersurfaces.
We consider closed characteristics $(\tau,y)$ on $\Sigma$, which are solutions of the following problem
$$
\left\{
\begin{aligned}
&\dot y=JN_\Sigma(y),\\
&y(\tau)=y(0).
\end{aligned}
\right.
$$
Here $J=\left(\begin{smallmatrix}0&-I_n\\ I_n&0\end{smallmatrix}\right)$, $I_n$ is the identity matrix in $\R^n$, $\tau>0$, and $N_\Sigma(y)$ is the outward normal vector of $\Sigma$ at $y$ normalized by the condition $N_\Sigma(y)\cdot y=1$.
Here $a\cdot b$ denotes the standard inner product of $a,b\in\R^{2n}$.
A closed characteristic $(\tau,y)$ is {\it prime} if $\tau$ is the minimal period of $y$.
Two closed characteristics $(\tau,y)$ and $(\sigma,z)$ are {\it geometrically distinct} if $y(\R)\neq z(\R)$.
We denote by $\mathcal T(\Sigma)$ the set of all geometrically distinct closed characteristics on $\Sigma$.
A closed characteristic $(\tau,y)$ is {\it non-degenerate} if $1$ is a Floquet multiplier of $y$ of precisely algebraic multiplicity $2$.
It is {\it hyperbolic} if $1$ is a double Floquet multiplier of $y$ and all the other Floquet multipliers are not on $\U=\{z\in\C\mid |z|=1\}$, i.e., the unit circle in the complex plane.
It is {\it elliptic} if all the Floquet multipliers of $y$ are on $\U$.
It is {\it irrationally elliptic} if it is elliptic and its linearized Poincar\'e map is suitably homotopic to the $\diamond$-product of $n-1$ rotation $(2\times 2)$ matrices with rotation angles being irrational multiples of $\pi$.
We call a hypersurface $\Sigma\in\mathcal{H}_{st}(2n)$ {\it non-degenerate} if all the closed characteristics on $\Sigma$, together with all of their iterations, are non-degenerate.

We call $\Sigma$ {\it symmetric} if $\Sigma=-\Sigma$.
We denote the set of all symmetric hypersurfaces in $\H_{st}(2n)$ (or $\H_{con}(2n)$) by $S\mathcal H_{\rm st}(2n)$ (or $S\mathcal H_{\rm con}(2n)$).
A closed characteristic $(\tau,y)$ on $\Sigma\in S\mathcal H_{\rm st}(2n)$ is called {\it symmetric} if $y(\R)=-y(\R)$.
If $(\tau,y)$ is prime and symmetric, then $y(t+\tau/2)=-y(t)$, $\forall\,t\in\R$.
We denote by $\mathcal T_s(\Sigma)$ the set of all geometrically distinct symmetric closed characteristics on $\Sigma$.

There is a long standing conjecture on the number of closed characteristics on compact star-shaped hypersurfaces in $\R^{2n}$:
\be \,^{\#}\T(\Sg)\ge n,\quad \forall\,\Sg\in\H_{st}(2n). \lb{conj}\ee
The global result $\,^{\#}\T(\Sigma)\ge1$ was proved by P. Rabinowitz \cite{Ra78} for $\Sigma\in\H_{st}(2n)$ and by A. Weinstein \cite{We78} for $\Sigma\in\H_{con}(2n)$ in 1978.
For further results on the multiplicity of closed characteristics on convex or star-shaped hypersurfaces, we refer to \cite{CGH16,EH87,EL87,HL02,LL16,WHL07} and the references therein.
In particular, Y. Long and C. Zhu \cite{LZ02} proved
$$
\,^{\#}\T(\Sigma)\ge \left[\frac{n}{2}\right]+1,\quad \forall\,\Sigma\in\H_{con}(2n),
$$
and obtained $\,^{\#}\T(\Sigma)\ge n$ under the non-degeneracy assumption.
More recently, E. Cineli, V. Ginzburg and B. G\"urel \cite{CGG24} proved the multiplicity conjecture \eqref{conj} in the dynamically convex setting, which includes the convex case.
D. Cristofaro-Gardiner, U. Hryniewicz, M. Hutchings and H. Liu \cite{CGHHL26} proved that $\,^{\#}\T(\Sigma)=2$ or $\infty$ for every smooth compact strictly star-shaped hypersurface $\Sigma\subset\R^4$.

Besides the multiplicity problem, the mean index plays an important role in the study of closed characteristics.
When there exist only finitely many geometrically distinct prime closed characteristics on $\Sigma \in \H_{st}(2n)$, the resonance identities of H. Liu, Y. Long and W. Wang \cite{LLW14} imply that at least one of them has positive mean index, while the others may have positive, zero or negative mean indices.
H. Duan, H. Liu, Y. Long and W. Wang \cite{DLLW24} generalized the common index jump theorem to the case where positive and negative mean indices may occur simultaneously, provided that all the mean indices are nonzero.

Motivated by analogous results for closed geodesics and by the examples of weakly non-resonant ellipsoids, H. Duan, H. Liu, Y. Long and W. Wang \cite{DLLW24} proposed the following conjecture \cite[Conjecture 1.5]{DLLW24}: 
if a compact star-shaped hypersurface carries only finitely many geometrically distinct closed characteristics, 
then none of them possesses zero mean index.
This conjecture is closely related to the multiplicity problem, because the nonzero mean index condition is the mean index assumption required in the applications of the generalized common index jump theorem in \cite{DLLW24}.
In Theorem 1.6 of \cite{DLLW24}, the above conjecture was proved for non-degenerate star-shaped hypersurfaces in $\R^6$.
The conjecture remains open in general dimensions.

Suppose now that $\Sigma\in S\mathcal H_{\rm st}(2n)$.
Symmetric closed characteristics on symmetric hypersurfaces have been studied by C. Liu, Y. Long and C. Zhu \cite{LLZ02}, and W. Wang \cite{Wa12}.
More recent results on symmetric Reeb orbits can be found, for example, in \cite{GM21,ALM24}.
In particular, H. Liu and Y. Long \cite{LL13,LL15} established resonance identities for symmetric closed characteristics in the convex and star-shaped settings.
They also proved that if a symmetric compact star-shaped hypersurface in $\R^4$ carries exactly two geometrically distinct closed characteristics, then both are elliptic; see Theorem 1.4 of \cite{LL15}.

Let $(\tau,y)$ be a prime symmetric closed characteristic and let $\gamma_y:[0,\tau/2]\to\Sp(2n)$ be the associated symplectic path used to define its symmetric indices in Section \ref{sec:2}.
Following H. Liu and Y. Long \cite{LL15}, we denote by
$$
\bar i(y)=i_{-1}(\gamma_y),\quad
\bar\nu(y)=\nu_{-1}(\gamma_y),\quad
\hat{\bar i}(y)=\lim_{m\to\infty}
\frac{\bar i(y^{2m-1})}{2m-1}
$$
the symmetric index, symmetric nullity and symmetric mean index, respectively.
For the definitions of the Maslov-type $(-1)$-index $i_{-1}$ and nullity $\nu_{-1}$, we refer to Chapter 5 of \cite{Lo02}, with $\omega=-1$.
Only odd iterates occur here because even iterates do not satisfy the symmetric boundary condition.
Moreover, by Lemma 5.5 of \cite{LL15},
$$
\hat{\bar i}(y)=\frac12\hat i(y),
$$
where $\hat i(y)$ denotes the ordinary mean index.

It is then natural to consider the symmetric counterpart of Conjecture 1.5 in \cite{DLLW24}: if a symmetric compact star-shaped hypersurface carries only finitely many geometrically distinct prime symmetric closed characteristics, 
that is, $\,^{\#}\mathcal T_s(\Sigma)<+\infty$, then $\hat{\bar i}(y)\neq0$ for every prime symmetric closed characteristic $y$.
The main result of this paper proves this statement in $\R^4$ under the non-degeneracy assumption in the symmetric sense on all odd iterates.

\begin{theorem}\label{thm:main}
Let $\Sigma=-\Sigma\subset\R^4$ be a $C^3$ compact strictly star-shaped symmetric hypersurface.
Assume that
\begin{enumerate}[label=(\roman*)]
\item $\,^{\#}\mathcal T_s(\Sigma)<+\infty$;
\item every odd iterate of every prime symmetric closed characteristic is non-degenerate in the symmetric sense, i.e., $\bar\nu(y^{2m-1})=1$ for every $m\in\N$.
\end{enumerate}
Then every prime symmetric closed characteristic $y$ satisfies
$$
\hat{\bar i}(y)\neq0.
$$
Consequently $\hat i(y)\neq0$ as well.
\end{theorem}

Theorem~\ref{thm:main} gives a symmetric counterpart in $\R^4$ of the result on the nonvanishing of mean indices in $\R^6$ proved by H. Duan, H. Liu, Y. Long and W. Wang \cite[Theorem~1.6]{DLLW24}.
The finiteness assumption here concerns only the prime symmetric closed characteristics, and non-symmetric closed characteristics may also exist on $\Sigma$.

\begin{remark}
The non-degeneracy assumption in Theorem \ref{thm:main} is weaker than the usual non-degeneracy assumption for closed characteristics.
Condition {\rm (ii)} only requires the odd iterates to be non-degenerate in the symmetric sense.
The nonzero critical points correspond only to odd iterates of symmetric closed characteristics.
Hence only these odd iterates are used in the index computation and in the local critical modules.
\end{remark}

Combining Theorem \ref{thm:main} with the resonance identities and Morse inequalities in Section~\ref{sec:2}, we obtain the following stability result.

\begin{theorem}\label{thm:two-symmetric-elliptic}
Let $\Sigma\in S\mathcal H_{\rm st}(4)$ satisfy $\,^{\#}\mathcal T_s(\Sigma)=2$.
Assume that all iterates of these two symmetric closed characteristics are non-degenerate.
Then both of them are irrationally elliptic.
\end{theorem}

Theorem~\ref{thm:two-symmetric-elliptic} extends the result of H. Liu and Y. Long \cite[Theorem~1.2]{LL13}, 
as well as that of D. Cristofaro-Gardiner, U. Hryniewicz, M. Hutchings and H. Liu \cite[Theorem~1.2]{CGHHL23},
to symmetric compact star-shaped hypersurfaces $\Sigma\subset\R^4$ with $\,^{\#}\mathcal T_s(\Sigma)=2$, 
provided that all iterates of the two symmetric closed characteristics are non-degenerate.

This paper is organized as follows.
In Section~\ref{sec:2}, we follow H. Liu and Y. Long \cite[Section~2]{LL15} to set up the variational structure for symmetric closed characteristics.
We also recall the resonance identities and Morse inequalities for symmetric closed characteristics.
In Section \ref{sec:3}, we give the symmetric index iteration formulae and prove a constant index property in dimension four.
In Section \ref{sec:4}, we adapt the arguments of Viterbo \cite[Section 7]{Vi89} to the symmetric setting and use these results to prove Theorem \ref{thm:main}.
Finally, we combine Theorem \ref{thm:main} with the resonance identities and Morse inequalities to prove Theorem \ref{thm:two-symmetric-elliptic}.

In this paper, let $\N$, $\N^0$, $\Z$, $\Q$, $\R$, and $\C$ denote the sets of natural integers, non-negative integers, integers, rational numbers, real numbers, and complex numbers respectively.
We define the functions $[a]=\max\{k\in\Z\mid k\le a\}$, $E(a)=\min\{k\in\Z\mid k\ge a\}$.
Denote by $a\cdot b$ and $|a|$ the standard inner product and norm in $\R^{2n}$.
For an $S^1$-space $X$, we denote by $X_{S^1}$ the homotopy quotient of $X$ by $S^1$, i.e., $X_{S^1}=S^\infty\times_{S^1}X$, where $S^\infty$ is the unit sphere in an infinite dimensional complex Hilbert space.
In this paper we use $\Q$ coefficients for all homological and cohomological modules.

\section{Critical point theory for symmetric closed characteristics}\label{sec:2}

In the following, we construct a variational structure of symmetric closed characteristics.
We first fix a $\Sigma\in S\mathcal H_{\rm st}(2n)$ and assume the following condition on $\cT_s(\Sigma)$:

\medskip
\noindent
$(\bf{F_s})$ There exist only finitely many geometrically distinct prime symmetric closed characteristics $\{(\tau_k,y_k)\}_{1\leq k\leq q}$ on $\Sigma$.
\medskip

Note that $(\tau,y)$ is some odd iteration of a prime symmetric closed characteristic if and only if it satisfies the problem
\begin{equation}\label{eq:symmetric-characteristic}
\left\{
\begin{aligned}
&\dot y(t)=JN_\Sigma(y(t)),\quad y(t)\in\Sigma,\quad \forall t\in\R,\\
&y(\tau/2)=-y(0).
\end{aligned}
\right.
\end{equation}
As in Definition 2.1 of \cite{LL15}, we introduce the following discrete subset of $\R^+$.
\begin{definition}\label{def:period-set}
Under the assumption $(\bf{F_s})$, the set of periods of symmetric closed characteristics on $\Sigma$ is defined by
$$
\per_s(\Sigma)=\{(2m-1)\tau_k\mid m\in\N,\ 1\leq k\leq q\}.
$$
\end{definition}
Note that in the above definition, the period set $\per_s(\Sigma)$ is defined only via odd iterations of $(\tau_k,y_k)$, because even iterations do not satisfy \eqref{eq:symmetric-characteristic} and do not yield any critical point of $\bar F_{a,K}$ on the space $\bar X$ below.

Now we construct the variational structure of symmetric closed characteristics as follows.
By Lemma 2.2 of \cite{LL15}, for any sufficiently small $\vartheta\in(0,1)$, there exists a function $\varphi\equiv\varphi_\vartheta\in C^\infty(\R,[0,+\infty))$ depending on $\vartheta$ which has $0$ as its unique critical point in $[0,+\infty)$ such that
\begin{equation}\label{eq:phi}
\varphi(0)=0=\varphi'(0),\quad
\varphi''(0)=1=\lim_{t\to0^+}\frac{\varphi'(t)}{t},
\end{equation}
and
\begin{equation}\label{eq:phi-monotone}
\frac{d}{dt}\left(\frac{\varphi'(t)}{t}\right)<0,\quad t>0,
\end{equation}
with $\lim_{t\to+\infty}\frac{\varphi'(t)}{t}<\vartheta$.
Moreover, we can choose $\alpha\in(1,2)$ sufficiently close to $2$ and $c\in(0,1)$ such that $\varphi(t)=ct^\alpha$ whenever $\varphi'(t)/t\in[\vartheta,1-\vartheta]$ and $t>0$.

Let $j:\R^{2n}\to\R$ be the gauge function of $\Sigma$, i.e., $j(\lambda x)=\lambda$, $x\in\Sigma$, $\lambda\geq0$ and $j(0)=0$.
Then $j\in C^3(\R^{2n}\setminus\{0\},\R)\cap C^0(\R^{2n},\R)$ and $\Sigma=j^{-1}(1)$.
Since $\Sigma=-\Sigma$, we also have $j(-x)=j(x)$, $x\in\R^{2n}$.
Let $\hat\tau=\inf \{s \mid s\in \per_s(\Sigma) \}$.
Let $T>0$ be fixed, and identify $S^1$ with $\R/T\Z$.
By Lemma 2.3 of \cite{LL15}, we have
\begin{lemma}\lb{tilde H}
Let $a>\hat\tau/T$, $\vartheta_a\in(0,\hat\tau/(aT))$ and $\varphi_a$ be a $C^\infty$ function associated to $\vartheta_a$ satisfying \eqref{eq:phi}-\eqref{eq:phi-monotone}.
Define the Hamiltonian function $\widetilde H_a(x)=a\varphi_a(j(x))$ and consider the fixed period system
\begin{equation}\label{eq:anti-Hamilton}
\left\{
\begin{aligned}
&\dot x(t)=J\widetilde H_a'(x(t)),\\
&x(T/2)=-x(0).
\end{aligned}
\right.
\end{equation}
Then the solutions of \eqref{eq:anti-Hamilton} are $x\equiv0$ and
\begin{equation}\label{eq:solutions}
x(t)=\rho y(\tau t/T),\quad
\frac{\varphi_a'(\rho)}{\rho}=\frac{\tau}{aT},
\end{equation}
where $(\tau,y)$ is a solution of \eqref{eq:symmetric-characteristic}.
In particular, nonzero solutions of \eqref{eq:anti-Hamilton} are in one-to-one correspondence with solutions of \eqref{eq:symmetric-characteristic} with period $\tau<aT$.
\end{lemma}

For technical reasons we want to further modify the Hamiltonian.
By Proposition 2.4 of \cite{LL15}, after choosing $\vartheta_a$ sufficiently small, we have
\begin{proposition}
For $a>\hat\tau/T$ and small $\varepsilon_a$, we choose small enough $\vartheta_a$ such that Lemma \ref{tilde H} holds.
Then there exists a function $\bar H_a$ on $\R^{2n}$ such that $\bar H_a(x)=\bar H_a(-x)$ for all $x\in\R^{2n}$, $\bar H_a$ is $C^1$ on $\R^{2n}$ and $C^3$ on $\R^{2n}\setminus\{0\}$, 
$\bar H_a=\widetilde H_a$ on $U_A\equiv\{x\mid \widetilde H_a(x)\leq A\}$ for some large $A$, and $\bar H_a(x)=\frac12\varepsilon_a|x|^2$ for $|x|$ large, where $\varepsilon_a>0$, 
and the solutions of the fixed period system
\begin{equation}\label{eq:modified-system}
\left\{
\begin{aligned}
&\dot x(t)=J\bar H_a'(x(t)),\\
&x(T/2)=-x(0),
\end{aligned}
\right.
\end{equation}
are the same as those of \eqref{eq:anti-Hamilton}.
\end{proposition}

As in Section 2 of \cite{LL15}, for any $a>\hat\tau/T$, we can choose some large constant $K=K(a)$ such that
$$
\bar H_{a,K}(x)=\bar H_a(x)+\frac12K|x|^2
$$
is a strictly convex function, that is,
\begin{equation}\label{eq:strict-convexity}
\bigl(\nabla\bar H_{a,K}(x)-\nabla\bar H_{a,K}(y),x-y\bigr)
\geq\frac{\epsilon}{2}|x-y|^2,
\end{equation}
for all $x,y\in\R^{2n}$ and some positive $\epsilon$.
Let $\bar H_{a,K}^*$ be the Fenchel dual of $\bar H_{a,K}$ defined by
$$
\bar H_{a,K}^*(y)=\sup\{x\cdot y-\bar H_{a,K}(x)\mid x\in\R^{2n}\}.
$$
Since $\bar H_a(x)=\bar H_a(-x)$ for all $x\in\R^{2n}$, $\bar H_{a,K}(x)=\bar H_{a,K}(-x)$ and $\bar H_{a,K}^*(x)=\bar H_{a,K}^*(-x)$ for all $x\in\R^{2n}$.

The dual action functional on
$$
\bar X=W_s^{1,2}(\R/T\Z,\R^{2n})
=\{x\in W^{1,2}(\R/T\Z,\R^{2n})\mid x(t+T/2)=-x(t),\ \forall t\in\R\}
$$
is defined by
\begin{equation}\label{eq:Fbar}
\bar F_{a,K}(x)=\int_0^{T/2}\left[
\frac12(J\dot x-Kx,x)+\bar H_{a,K}^*(-J\dot x+Kx)
\right]dt.
\end{equation}
Then $\bar F_{a,K}\in C^{1,1}(\bar X,\R)$, but $\bar F_{a,K}$ is not $C^2$.

Assume
\begin{equation}\label{eq:K-nonresonance}
\frac{KT}{2\pi}\notin2\Z-1.
\end{equation}
Then, as in Lemma 2.5 of \cite{LL15}, the map $x\longmapsto-J\dot x+Kx$ is a Hilbert space isomorphism between $\bar X=W_s^{1,2}(\R/T\Z,\R^{2n})$ and
$$
\bar E=L_s^2(\R/T\Z,\R^{2n})
=\{u\in L^2(\R/T\Z,\R^{2n})\mid u(t+T/2)=-u(t),\ \text{a.e. }t\in\R\}.
$$
Consequently, by Lemma 2.5 of \cite{LL15}, we have the following lemma.
\begin{lemma}\lb{LL15 2.5}
Assume $\frac{KT}{2\pi}\notin2\Z-1$.
Then $x$ is a critical point of $\bar F_{a,K}$ if and only if it is a solution of \eqref{eq:modified-system}.
\end{lemma}
From Lemma \ref{LL15 2.5}, we know that the critical points of $\bar F_{a,K}$ are independent of $K$.
By Proposition 2.6 of \cite{LL15}, for every nonzero critical point $x_a$ of $\bar F_{a,K}$, its critical value satisfies $\bar F_{a,K}(x_a)<0$ and is independent of $K$.
More precisely, by \eqref{eq:solutions}, we have $x_a=\rho_a y(\tau t/T)$ with $\varphi_a'(\rho_a)/\rho_a = \tau/(aT)$, and
$$
\bar F_{a,K}(x_a)=\frac14a\varphi_a'(\rho_a)\rho_aT-\frac12a\varphi_a(\rho_a)T<0.
$$

We know that when \eqref{eq:K-nonresonance} holds, the map $x\mapsto-J\dot x+Kx$ is a Hilbert space isomorphism between $\bar X$ and $\bar E$.
We denote its inverse by $\bar M_K$ and define the functional
\begin{equation}\label{eq:Psibar}
\bar\Psi_{a,K}(u)=\int_0^{T/2}\left[
-\frac12(\bar M_Ku,u)+\bar H_{a,K}^*(u)
\right]dt,\quad \forall u\in\bar E.
\end{equation}
Then $x\in\bar X$ is a critical point of $\bar F_{a,K}$ if and only if $u=-J\dot x+Kx$ is a critical point of $\bar\Psi_{a,K}$.

We have a natural $S^1$-action on $\bar X$ or $\bar E$ defined by
$$
\theta\cdot x(t)=x(\theta+t),\quad \forall\theta\in S^1,\ t\in\R.
$$
Then the functionals $\bar F_{a,K}$ and $\bar\Psi_{a,K}$ are $S^1$-invariant and are even since $\bar H_{a,K}^*$ is even.
If a nonzero critical point corresponds to the $m$-th iteration of a prime symmetric closed characteristic, where $m$ is odd, then its isotropy group is $\Z_m$ and its critical orbit is isomorphic to $S^1/\Z_m$.
For any $\kappa\in\R$, set
$$
\bar\Lambda_{a,K}^{\kappa}=\{u\in\bar E\mid \bar\Psi_{a,K}(u)\leq\kappa\},\quad
\bar X_{a,K}^{\kappa}=\{x\in\bar X\mid \bar F_{a,K}(x)\leq\kappa\}.
$$
These sets are $S^1$-invariant. 
When $a$ and $K$ are fixed, we write $\bar X^\kappa$ for $\bar X_{a,K}^\kappa$.

As in Definition 2.7 of \cite{LL15}, we have
\begin{definition}
Suppose $u$ is a nonzero critical point of $\bar\Psi_{a,K}$.
Then the formal Hessian of $\bar\Psi_{a,K}$ at $u$ is defined by
\begin{equation}\label{eq:Hessian}
\bar Q_{a,K}(v)=\int_0^{T/2}
\left(-\bar M_Kv\cdot v+\bar H_{a,K}^{*\prime\prime}(u)v\cdot v\right)dt,
\end{equation}
which defines an orthogonal splitting $\bar E=\bar E^-\oplus\bar E^0\oplus\bar E^+$ of $\bar E$ into negative, zero and positive subspaces.
The index and nullity of $u$ are defined by
$$
\bar i_K(u)=\dim\bar E^-,\quad
\bar\nu_K(u)=\dim\bar E^0,
$$
respectively.
\end{definition}

Similarly, we define the index and nullity of $x=\bar M_Ku$ for $\bar F_{a,K}$ respectively by
$$
\bar i_K(x)=\bar i_K(u),\quad \bar\nu_K(x)=\bar\nu_K(u).
$$

By the argument used before Theorem 2.8 of \cite{LL15}, we have that $v\in\bar E$ belongs to the null space of $\bar Q_{a,K}$ if and only if $z=\bar M_Kv$ is a solution of the linearized system
$$
\left\{
\begin{aligned}
&\dot z(t)=J\bar H_a''(x(t))z(t),\\
&z(T/2)=-z(0).
\end{aligned}
\right.
$$
Thus the nullity in \eqref{eq:Hessian} is independent of $K$, which we denote by $\bar{\nu}(x)\equiv \bar\nu_K(u)=\bar\nu_K(x)$.

Consider the linear Hamiltonian system
\begin{equation}\label{eq:linear-Hamiltonian}
  \left\{
\begin{aligned}
&\dot\xi(t)=JA_a(t)\xi(t),\\
&A_a(t+T/2)=A_a(t),\quad t\in\R,
\end{aligned}
\right.
\end{equation}
where $A_a(t)=\bar H_a''(x(t))$.
The second relation in \eqref{eq:linear-Hamiltonian} follows from the evenness of $\bar H_a$ and the identity $x(t+T/2)=-x(t)$.
Let $\gamma_x:[0,T/2]\to\Sp(2n)$ with $\gamma_x(0)=I_{2n}$ be the fundamental solution of \eqref{eq:linear-Hamiltonian}.
In the following we set
$$
e(K)=2n\left(\left[\frac{KT}{2\pi}\right]-\left[\frac{KT}{4\pi}\right]\right).
$$
By Theorem 2.8 of \cite{LL15}, we have
\begin{equation}\label{eq:Morse-Maslov}
\bar i_K(u)=\bar i_K(x)=e(K)+i_{-1}(\gamma_x),\quad
\bar\nu_K(u)=\bar\nu_K(x)=\nu_{-1}(\gamma_x),
\end{equation}
where $\nu_{-1}(\gamma_x)=\dim\ker(\gamma_x(T/2)+I)$.
The index $i_{-1}(\gamma_x)$ and nullity $\nu_{-1}(\gamma_x)$ are independent of $K$.

In the following, we give the definitions of the index and nullity for symmetric closed characteristics.
For every solution $(\tau,y)$ of \eqref{eq:symmetric-characteristic}, that is, every odd iterate of a prime symmetric closed characteristic, let $aT>\tau$ and choose $\varphi_a$ to satisfy \eqref{eq:phi}-\eqref{eq:phi-monotone}.
Determine $\rho$ uniquely by $\varphi_a'(\rho)/\rho=\tau/(aT)$, and let $x(t)=\rho y(\tau t/T)$.
Then we define the index $\bar i(\tau,y)$ and nullity $\bar\nu(\tau,y)$ of $(\tau,y)$ by
$$
\bar i(\tau,y)=i_{-1}(\gamma_x),\quad \bar\nu(\tau,y)=\nu_{-1}(\gamma_x).
$$
By Proposition 2.12 of \cite{LL15}, the index and nullity are well defined and are independent of the choice of $aT>\tau$ and $\varphi_a$ satisfying \eqref{eq:phi}-\eqref{eq:phi-monotone}.
Moreover, when $\bar H_a$ is $\alpha$-homogeneous near the image set of $x$, they coincide with those defined for the Hamiltonian $H(x)=j(x)^\alpha$ for all $x\in\R^{2n}$.
In particular, $1\leq\bar\nu(\tau,y)\leq2n-1$.
For a prime symmetric closed characteristic $(\tau,y)$ and $m\in\N$, we write $y^{2m-1}\equiv((2m-1)\tau,y)$ and $\bar i(y^{2m-1})=\bar i((2m-1)\tau,y)$, $\bar\nu(y^{2m-1})=\bar\nu((2m-1)\tau,y)$.

By Proposition 2.13 of \cite{LL15}, we have
\begin{proposition}
$\bar\Psi_{a,K}$ satisfies the Palais-Smale condition on $\bar E$, and $\bar F_{a,K}$ satisfies the Palais-Smale condition on $\bar X$ when \eqref{eq:K-nonresonance} holds.
\end{proposition}

We shall also use the finite-dimensional reduction in \cite[Lemma 2.11]{LL15}.
For any $\hat\tau/T<a_1<a_2<+\infty$, fix $K$ such that $\bar H_{a,K}$ satisfies \eqref{eq:strict-convexity} with the same $\epsilon>0$ for every $a\in[a_1,a_2]$.
Assume that $\{\bar\Psi_{a,K}\}_{a\in[a_1,a_2]}$ is a continuous family of functionals in the sense of Remark~2.6 of \cite{LLW14}.
Then there exist a finite-dimensional $S^1$-invariant subspace $\bar G$ of $\bar E$ and a family of $S^1$-equivariant maps $\bar h_a:\bar G\to\bar G^\perp$ such that, for $g\in\bar G$, each function
$$
h\longmapsto\bar\Psi_{a,K}(g+h)
$$
has $\bar h_a(g)$ as the unique minimizer in $\bar G^\perp$.
Let
$$
\bar\gamma_{a,K}(g)=\bar\Psi_{a,K}(g+\bar h_a(g)).
$$
Then $\bar\gamma_{a,K}$ is $C^1$ and $S^1$-invariant on $\bar G$, and $g_a$ is a critical point of $\bar\gamma_{a,K}$ if and only if $g_a+\bar h_a(g_a)$ is a critical point of $\bar\Psi_{a,K}$.
If $g_a$ is a nonzero critical point of $\bar\gamma_{a,K}$, then $\bar\gamma_{a,K}$ is $C^2$ in a neighborhood of the critical orbit $S^1\cdot g_a$, and the index and nullity coincide with those of $\bar\Psi_{a,K}$.
Moreover, the natural embedding of the reduced sublevel into the corresponding sublevel of $\bar\Psi_{a,K}$ is an $S^1$-equivariant homotopy equivalence.
The map $a\mapsto\bar\gamma_{a,K}$ is continuous in the $C^1$ topology, and its Hessian depends continuously on $a$ in a neighborhood of a nonzero critical orbit.

By Lemma 2.9 and Corollary 2.10 of \cite{LL15}, for $b>0$ sufficiently large, if $\bar X^{-\infty}=\bar X_{a,K}^{-b}$, then the global equivariant cohomological structure of $(\bar X,\bar X^{-\infty})$ is given by
$$
H_{S^1}^q(\bar X,\bar X^{-\infty};\Q) \cong H^{q-e(K)}(\mathbf CP^\infty;\Q),
$$
which is one-dimensional for all even $q\geq e(K)$ and zero otherwise.

Now for a critical point $u$ of $\bar\Psi_{a,K}$ and the corresponding $x=\bar M_Ku$ of $\bar F_{a,K}$, let
$$
\bar\Lambda_{a,K}(u)=\{w\in\bar E\mid \bar\Psi_{a,K}(w)\leq\bar\Psi_{a,K}(u)\}.
$$
Denote by $\operatorname{crit}(\bar\Psi_{a,K})$ the set of critical points of $\bar\Psi_{a,K}$.
Because $\bar\Psi_{a,K}$ is $S^1$-invariant, $S^1\cdot u$ is a critical orbit if $u\in\operatorname{crit}(\bar\Psi_{a,K})$.
Under the condition $(F_s)$, the number of critical orbits of $\bar\Psi_{a,K}$ is finite.
Hence as usual we can make the following definition.

\begin{definition}\label{def:critical-module}
Suppose $u$ is a nonzero critical point of $\bar\Psi_{a,K}$, and $\cN$ is an $S^1$-invariant open neighborhood of $S^1\cdot u$ such that $\operatorname{crit}(\bar\Psi_{a,K})\cap(\bar\Lambda_{a,K}(u)\cap\cN)=S^1\cdot u.$
Then the $S^1$-critical modules of $S^1\cdot u$ are defined by
$$
C_{S^1,q}(\bar\Psi_{a,K},S^1\cdot u)
=H_q\left((\bar\Lambda_{a,K}(u)\cap\cN)_{S^1},
((\bar\Lambda_{a,K}(u)\setminus S^1\cdot u)\cap\cN)_{S^1};\Q\right).
$$
\end{definition}
Similarly, we define the $S^1$-critical modules $C_{S^1,q}(\bar F_{a,K},S^1\cdot x)$ of $S^1\cdot x$ for $\bar F_{a,K}$.
By \eqref{eq:Fbar} and \eqref{eq:Psibar}, we have
$$
C_{S^1,q}(\bar\Psi_{a,K},S^1\cdot u) \cong C_{S^1,q}(\bar F_{a,K},S^1\cdot x).
$$

By Proposition 2.15 of \cite{LL15}, the critical module $C_{S^1,q}(\bar F_{a,K},S^1\cdot x)$ is independent of the choices of $a$, $\varphi_a$ and $\bar H_a$ when the critical points correspond to the same symmetric closed characteristic.
By Theorem 2.16 of \cite{LL15}, the $S^1$-critical module $C_{S^1,e(K)+l}(\bar F_{a,K},S^1\cdot x)$ is independent of the choice of $K$, i.e.,
$$
C_{S^1,e(K)+l}(\bar F_{a,K},S^1\cdot x)\cong
C_{S^1,e(K')+l}(\bar F_{a,K'},S^1\cdot x), \quad l\in\Z,
$$
for $K$ and $K'$ satisfying the above conditions.

Finally, let $m\in\N$, choose $a$ such that $aT>(2m-1)\tau$, and let $u^{2m-1}$ be the critical point of $\bar\Psi_{a,K}$ corresponding to the odd iteration $y^{2m-1}$.
If $y^{2m-1}$ is non-degenerate in the symmetric sense, i.e., $\bar\nu(y^{2m-1})=1$, then by Proposition 3.1 of \cite{LL15}, together with \eqref{eq:Morse-Maslov}, Proposition 2.15 and Theorem 2.16 of \cite{LL15}, we have
\begin{equation}\label{eq:local-module}
C_{S^1,q}(\bar\Psi_{a,K},S^1\cdot u^{2m-1})\cong
\begin{cases}
\Q,&q=e(K)+\bar i(y^{2m-1}),\\
0,&\text{otherwise}.
\end{cases}
\end{equation}

We have the following mean index identities for symmetric closed characteristics.

\begin{theorem}[cf. Theorem 1.1 of {\cite{LL15}}]\label{thm:symmetric-resonance}
Suppose that $\Sigma\in S\mathcal H_{\rm st}(2n)$ satisfies $\,^{\#}\mathcal T_s(\Sigma)<+\infty$.
Denote by $\{(\tau_j,y_j)\}_{1\leq j\leq k}$ all the geometrically distinct prime symmetric closed characteristics.
Then the following identities hold
\begin{equation}\label{eq:symmetric-resonance}
 \sum_{\substack{1\leq j\leq k\\ \hat{\bar i}(y_j)>0}}
 \frac{\hat{\bar\chi}(y_j)}{\hat{\bar i}(y_j)}=1,\quad
 \sum_{\substack{1\leq j\leq k\\ \hat{\bar i}(y_j)<0}}
 \frac{\hat{\bar\chi}(y_j)}{\hat{\bar i}(y_j)}=0,
\end{equation}
where $\hat{\bar\chi}(y)\in\Q$ is the average Euler characteristic of $y$ given by Definition 3.7 and Remark 3.8 of \cite{LL15}.
\end{theorem}

In particular, if $\bar\nu(y^{2m-1})=1$ for every $m\in\N$, then (3.20) of \cite{LL15} gives
\begin{equation}\label{eq:symmetric-average-euler}
 \hat{\bar\chi}(y)=(-1)^{\bar i(y)}.
\end{equation}

Let $\bar F_{a,K}$ be the functional defined in \eqref{eq:Fbar} for sufficiently large $a$ and $K$ satisfying $aT\notin\per_s(\Sigma)$ and \eqref{eq:K-nonresonance}.
Choose $\epsilon>0$ small enough that $[-\epsilon,0)$ contains no critical values of $\bar F_{a,K}$.
For $b$ large enough, the normalized Morse series of $\bar F_{a,K}$ in $\bar X^{-\epsilon}\setminus\bar X^{-b}$ is defined, as usual, by
\begin{equation}\label{eq:normalized-symmetric-morse-series}
 \bar M_a(t)=\sum_{\substack{q\geq0, 1\leq j\leq p}}
 \dim C_{S^1,q}(\bar F_{a,K},S^1\cdot v_j)t^{q-e(K)},
\end{equation}
where we denote by $\{S^1\cdot v_1,\ldots,S^1\cdot v_p\}$ the critical orbits of $\bar F_{a,K}$ with critical values less than $-\epsilon$.

The Poincar\'e series of $H_{S^1,*}(\bar X,\bar X^{-\epsilon};\Q)$ is $t^{e(K)}\bar Q_a(t)$.
According to Theorem~4.1 and (4.4) of \cite{LL15}, if we set $\bar Q_a(t)=\sum_{k\in\Z}q_kt^k$, then
$$
 q_k=0,\quad \forall k\in\mathring I,
$$
where $I$ is an interval of $\Z$ such that
$I\cap[\bar i(\tau,y),\bar i(\tau,y)+\bar\nu(\tau,y)-1]=\emptyset$
for all symmetric closed characteristics $(\tau,y)$ on $\Sigma$ satisfying \eqref{eq:symmetric-characteristic} with $\tau\geq aT$.
Here $\mathring I$ denotes $I$ with its endpoints removed.
Then by (4.6) of \cite{LL15}, we have
$$
 \bar M_a(t)-\frac{1}{1-t^2}+\bar Q_a(t)=(1+t)\bar U_a(t),
$$
where $\bar U_a(t)=\sum_{i\in\Z}u_{i}t^i$ is a Laurent series with nonnegative coefficients.
If there is no symmetric closed characteristic $(\tau,y)$ with $\hat{\bar i}(y)=0$, then by (4.18) of \cite{LL15} we have
\begin{equation}\label{eq:symmetric-morse-series}
 \bar M(t)-\frac{1}{1-t^2}=(1+t)\bar U(t),
\end{equation}
where $\bar M(t)=\sum_{p\in\Z}M_pt^p$ denotes $\bar M_a(t)$ as $a$ tends to infinity, and $\bar U(t)=\sum_{p\in\Z}u_pt^p$ denotes $\bar U_a(t)$ as $a$ tends to infinity and possesses only non-negative coefficients.
In addition, we also denote by $\beta_p$ the coefficient of $t^p$ of $\frac{1}{1-t^2}=\sum_{p\in\Z}\beta_pt^p$, i.e. there holds $\beta_p=1$ for $p\in2\N^0$ and $\beta_p=0$ for $p\notin2\N^0$.
Comparing coefficients in \eqref{eq:symmetric-morse-series} gives
\begin{equation}\label{eq:symmetric-morse-coefficients}
 M_p-\beta_p=u_p+u_{p-1}\geq0,\quad p\in\Z.
\end{equation}

For any positive integer $n_1$ and any nonnegative integer $n_2$, it follows from \eqref{eq:symmetric-morse-series} that
$$
 \begin{aligned}
 \sum_{p=-2n_1+1}^{2n_2+1}M_pt^p-\sum_{p=-2n_1+1}^{2n_2+1}\beta_pt^p
 &=(1+t)\sum_{p=-2n_1}^{2n_2+1}u_pt^p\\
 &\quad-u_{2n_2+1}t^{2n_2+2}-u_{-2n_1}t^{-2n_1},
 \end{aligned}
$$
which, through letting $t=-1$, yields the following Morse inequality
\begin{equation}\label{eq:symmetric-morse-odd}
 \sum_{p=-2n_1+1}^{2n_2+1}(-1)^pM_p \leq\sum_{p=-2n_1+1}^{2n_2+1}(-1)^p\beta_p.
\end{equation}
Similarly we have
\begin{equation}\label{eq:symmetric-morse-even}
 \sum_{p=-2n_1}^{2n_2}(-1)^pM_p \geq\sum_{p=-2n_1}^{2n_2}(-1)^p\beta_p.
\end{equation}
If $M_p=0$ for every $p<0$, then \eqref{eq:symmetric-morse-odd} and \eqref{eq:symmetric-morse-even}, together with $\beta_p=0$ for $p<0$, give
\begin{equation}\label{eq:stability-Morse}
 \sum_{p=0}^{k}(-1)^{k-p}M_p \geq \sum_{p=0}^{k}(-1)^{k-p}\beta_p, \quad k\in\N^0.
\end{equation}

\section{Symmetric index iteration}\label{sec:3}

In this section, we study the index iteration of symmetric closed characteristics.
The Bott iteration formula and the precise iteration formulae for the Maslov-type index were established by Y. Long \cite{Lo99,Lo00}. See also Chapters 8-9 of \cite{Lo02}.

We now consider the symmetric index iteration in $\R^4$.
Let $(\tau,y)$ be a prime symmetric closed characteristic on $\Sigma$.
As in Section \ref{sec:2}, we compute its symmetric index and nullity using $H=j^\alpha$, $1<\alpha<2$.
Let $x$ be the corresponding $\tau$-periodic Hamiltonian orbit, and let $\gamma:[0,\tau/2]\to\Sp(4)$, $\gamma(0)=I_4$, be its associated symplectic path.
As in (3.9) of \cite{DLLW24}, for $m\in\N$ we define the iterated path $\gamma^m:[0,m\tau/2]\to\Sp(4)$ by
$$
 \gamma^m(t)=\gamma(t-j\tau/2)[\gamma(\tau/2)]^j,
 \quad \frac{j\tau}{2}\leq t\leq\frac{(j+1)\tau}{2},
 \quad j=0,\ldots,m-1.
$$
The odd iterate $y^{2m-1}$ corresponds to $\gamma^{2m-1}$.
By Theorem 2.8 and Definition 3.6 of \cite{LL15}, the symmetric index and nullity of $y^{2m-1}$ are the Maslov-type $(-1)$-index and nullity of the associated symplectic path $\gamma^{2m-1}$.
Therefore, for every $m\in\N$, we have
\begin{equation}\label{eq:symmetric-index-nullity}
\begin{aligned}
 \bar i(y^{2m-1})=i_{-1}(\gamma^{2m-1}),\quad
 \bar\nu(y^{2m-1})=\nu_{-1}(\gamma^{2m-1})=\dim\ker\bigl([\gamma(\tau/2)]^{2m-1}+I_4\bigr),
\end{aligned}
\end{equation}
where the last equality follows from the definition of the $(-1)$-nullity.
The Bott formula gives
$$
 i_1(\gamma^{4m-2})=i_1(\gamma^{2m-1})+i_{-1}(\gamma^{2m-1}), \quad m\in\N,
$$
where $i_1$ denotes the ordinary Maslov-type index defined in Chapter 5 of \cite{Lo02}.
Consequently,
\begin{equation}\label{eq:advisor-difference}
 \bar i(y^{2m-1})=i_1(\gamma^{4m-2})-i_1(\gamma^{2m-1}), \quad m\in\N.
\end{equation}

We use the following basic normal forms.
As in \cite{Lo00}, denote by
$$
N_1(\lambda,b)=
 \begin{pmatrix}
 \lambda&b\\
 0&\lambda
 \end{pmatrix}, \quad \text{for }\lambda=\pm 1, \quad b\in \R.
$$
$$
D(\lambda)=
\begin{pmatrix}
\lambda & 0\\
0 & \lambda^{-1}
\end{pmatrix}, \quad \text{for }\lambda\in\R\setminus\{0,\pm1\}.
$$

$$
R(\theta)=
\begin{pmatrix}
\cos\theta & -\sin\theta\\
\sin\theta & \cos\theta
\end{pmatrix},\quad \text{for }\theta\in(0,\pi)\cup(\pi,2\pi).
$$

As in \cite{Lo00}, the $\diamond$-sum  of any two real matrices is defined by
$$
\begin{pmatrix}
A_1 & B_1\\
C_1 & D_1
\end{pmatrix}_{2i\times2i}
\diamond
\begin{pmatrix}
A_2 & B_2\\
C_2 & D_2
\end{pmatrix}_{2j\times2j}
=
\begin{pmatrix}
A_1 & 0   & B_1 & 0\\
0   & A_2 & 0   & B_2\\
C_1 & 0   & D_1 & 0\\
0   & C_2 & 0   & D_2
\end{pmatrix}.
$$
For every $M\in\Sp(2n)$, the homotopy set $\Omega(M)$ of $M$ in $\Sp(2n)$ is defined by
$$
\Omega(M)=
\left\{
N\in\Sp(2n)
\,\middle|\,
\begin{array}{l}
\sigma(N)\cap\U=\sigma(M)\cap\U\equiv\Gamma,  \nu_\omega(N)=\nu_\omega(M), \forall\,\omega\in\Gamma
\end{array}
\right\},
$$
where $\sigma(M)$ denotes the spectrum of $M$, and $ \nu_\omega(M) \equiv \dim_{\C}\ker_{\C}(M-\omega I)$ for $\omega\in\U$.
The component $\Omega^0(M)$ of $M$ in $\Sp(2n)$ is defined as the path-connected component of $\Omega(M)$ containing $M$.

In the symmetric case considered here, we need the corresponding decomposition of $\gamma(\tau/2)$.
Following the same idea as in Lemma 15.2.3 and Lemma 15.2.4 of \cite{Lo02}, see also (6.12) of \cite{LL15}, we obtain the following half-period version for symmetric closed characteristics.
\begin{lemma}\label{lem:first-block}
There exist $S\in\Sp(4)$ and $M\in\Sp(2)$ such that
\begin{equation}\label{eq:half-period-decomposition}
\gamma(\tau/2)=S^{-1}\bigl(N_1(-1,-1)\diamond M\bigr)S.
\end{equation}
\end{lemma}

We next give the precise odd-iteration formula in $\R^4$.
This is the symmetric counterpart of (5.3) of \cite{DLLW24}.
Assume now that
\begin{equation}\label{eq:odd-nondegenerate}
 \bar\nu(y^{2m-1})=1,\quad m\in\N.
\end{equation}
By \eqref{eq:symmetric-index-nullity} and \eqref{eq:half-period-decomposition}, for every $m\in\N$ we have
$$
N_1(-1,-1)^{2m-1}=N_1(-1,-(2m-1)),
$$
and hence
$$
 \bar\nu(y^{2m-1})=1+\dim\ker(M^{2m-1}+I_2).
$$
It follows from \eqref{eq:odd-nondegenerate} that
\begin{equation}\label{eq:ker-M}
 \ker(M^{2m-1}+I_2)=\{0\},\quad m\in\N.
\end{equation}

Since $M$ has no eigenvalue $-1$, the normal form of $\gamma(\tau/2)$ contains one block $N_1(-1,-1)$ and no blocks $N_1(-1,1)$ or $-I_2$. Thus, in the notation of Theorem 1.3 of \cite{Lo00}, $q_+=1$ and $q_-=q_0=0$.
By Theorem 1.8.10 of \cite{Lo02}, $M$ can be connected within $\Omega^0(M)$ to a diamond product of basic normal forms as
\begin{equation}\label{eq:M-normal-form}
\begin{aligned}
 &N_1(1,1)^{\diamond p_-}
 \diamond I_2^{\diamond p_0}
 \diamond N_1(1,-1)^{\diamond p_+}
 \diamond R(\theta_1)\diamond\cdots\diamond R(\theta_r)
 \diamond D(\pm2)^{\diamond s},
\end{aligned}
\end{equation}
where $p_-,p_0,p_+,r,s\in\N^0$, $\theta_j\in(0,\pi)\cup(\pi,2\pi)$ for $1\leq j\leq r$, and $D(\pm2)$ denotes either $D(2)$ or $D(-2)$.
Since $M\in\Sp(2)$, the dimensions of the blocks in \eqref{eq:M-normal-form} satisfy
\begin{equation}\label{eq:M-dimension-relation}
 p_-+p_0+p_++r+s=1.
\end{equation}
Moreover, condition \eqref{eq:odd-nondegenerate} implies
\begin{equation}\label{eq:angle-nonresonance}
 \frac{(2m-1)\theta_j}{2\pi}\notin\Z+\frac12, \quad m\in\N,\quad 1\le j\le r.
\end{equation}
Indeed, if the number in \eqref{eq:angle-nonresonance} belonged to $\Z+\frac12$, then $-1$ would be an eigenvalue of $M^{2m-1}$.
This contradicts \eqref{eq:ker-M}.

We first write the ordinary iteration formula for the half-period path $\gamma$.
Note that Theorem 8.3.1 of \cite{Lo02} gives
\begin{equation}\label{eq:full-iteration}
\begin{aligned}
 i_1(\gamma^k)={}&k\bigl(i_1(\gamma)+p_-+p_0-r\bigr)
 +2\sum_{j=1}^{r}E\left(\frac{k\theta_j}{2\pi}\right)-r-p_--p_0-\frac{1+(-1)^k}{2}, \quad k\in\N.
\end{aligned}
\end{equation}
If $x\notin\Z+\frac12$, then $E(2x)-E(x)=[2x]-[x]$.
Combining \eqref{eq:full-iteration} with \eqref{eq:advisor-difference}, \eqref{eq:angle-nonresonance} and the above identity, we obtain the precise iteration formula
\begin{equation}\label{eq:R4-precise-odd-formula}
\begin{aligned}
 \bar i(y^{2m-1}) =&i_1(\gamma^{4m-2})-i_1(\gamma^{2m-1}) \\
 ={}&(2m-1)\bigl(i_1(\gamma)+p_-+p_0-r\bigr)-1 \\
 &+2\sum_{j=1}^{r}
 \left( \left[\frac{(4m-2)\theta_j}{2\pi}\right] -\left[\frac{(2m-1)\theta_j}{2\pi}\right] \right), \quad m\in\N.
\end{aligned}
\end{equation}
By \eqref{eq:R4-precise-odd-formula}, we have
\begin{equation}\label{eq:symmetric-mean-formula}
 \hat{\bar i}(y)=\lim_{m\to+\infty}\frac{\bar i(y^{2m-1})}{2m-1}=i_1(\gamma)+p_-+p_0-r+\sum_{j=1}^{r}\frac{\theta_j}{\pi}.
\end{equation}

We now consider the case $\hat{\bar i}(y)=0$.

\begin{lemma}\label{lem:R4-zero-mean-classification}
Let $y$ be a prime symmetric closed characteristic in $\R^4$.
Assume that \eqref{eq:odd-nondegenerate} holds.
If $\hat{\bar i}(y)=0$, then $r=0$ and
$$
 \bar i(y^{2m-1})=-1, \quad m\in\N.
$$
\end{lemma}

\begin{proof}
By \eqref{eq:symmetric-mean-formula}, the zero mean index condition is equivalent to
\begin{equation}\label{eq:zero-mean-equation}
 i_1(\gamma)+p_-+p_0-r +\sum_{j=1}^{r}\frac{\theta_j}{\pi}=0.
\end{equation}
The dimension relation \eqref{eq:M-dimension-relation} gives $0\le r\le1$.

Suppose first that $r=1$.
In this case, \eqref{eq:zero-mean-equation} gives
$$
 \frac{\theta_1}{\pi} =1-i_1(\gamma)-p_--p_0\in\Z.
$$
This is impossible because $\theta_1\in(0,\pi)\cup(\pi,2\pi)$.
Therefore $r\ne1$, and hence $r=0$.

Now \eqref{eq:zero-mean-equation} gives
$$
 i_1(\gamma)+p_-+p_0=0.
$$
Together with \eqref{eq:R4-precise-odd-formula}, we obtain
$$
 \bar i(y^{2m-1})=-1, \quad m\in\N.
$$
This completes the proof.
\end{proof}

\section{Proofs of the main results}\label{sec:4}

We first prove Theorem \ref{thm:main} by studying the contribution of the origin to the functional $\bar F_{a,K}$ and adapting the arguments of Viterbo \cite[Section 7]{Vi89} to the symmetric setting.
We then prove Theorem \ref{thm:two-symmetric-elliptic} using the resonance identities and Morse inequalities.

As proved at the beginning of Section 7 of \cite{Vi89}, the contribution of the origin depends only on the Hamiltonian in a neighborhood of the origin.
Thus, for the proof of Theorem \ref{thm:main}, we use Hamiltonian functions that are homogeneous of degree two, i.e.,
$$
 \bar H_a(x)=\frac a2j(x)^2,\quad \bar H_{a,K}(x)=\bar H_a(x)+\frac K2|x|^2.
$$
Until the end of the proof of Theorem \ref{thm:main}, the notation from Section \ref{sec:2} refers to this model.
We set
$$
 \bar F_{a,K}(x)=\int_0^{T/2} \left[\frac12(J\dot x-Kx,x) +(\bar H_{a,K})^*(-J\dot x+Kx)\right]ds.
$$
For $c\in\R$, denote
$$
 \bar X^c(a,K)=\{x\in\bar X:\bar F_{a,K}(x)\leq c\},\quad
 \bar X^-(a,K)=\{x\in\bar X:\bar F_{a,K}(x)<0\}.
$$
Since by assumption $(\bf{F_s})$ the set of periods $\operatorname{per}_s(\Sigma)$ is discrete, we can choose $a>0$ such that $aT\notin\operatorname{per}_s(\Sigma)$.
We can also choose $K$ sufficiently large such that $\bar H_{a,K}$ is strictly convex and $KT/(2\pi)\notin2\Z-1$.

A nonzero solution of
$$
 \dot x=J\nabla\bar H_a(x),\quad x(T/2)=-x(0),
$$
has the form $x(t)=\rho y(at)$ where $y$ is a solution of \eqref{eq:symmetric-characteristic} with period $aT$.
But our choice of $a$ excludes these solutions.
Since $-J\frac d{ds}+K$ is invertible, the origin is the only critical point of $\bar F_{a,K}$.
For sufficiently small $\epsilon>0$, the group
$$
 H_{S^1,*}(\bar X^\epsilon(a,K),\bar X^{-\epsilon}(a,K);\Q)
$$
is the contribution of the origin for $\bar F_{a,K}$.
The following theorem is the symmetric version of Theorem 7.1 of \cite{Vi89}; see also Theorem 4.1 of \cite{LL15}. 
\begin{theorem}\label{thm:symmetric-Viterbo}
Let $a>0$ satisfy $aT\notin\operatorname{per}_s(\Sigma)$, and let $K$ be a sufficiently large constant satisfying $\frac{KT}{2\pi}\notin2\Z-1$ such that $\bar H_{a,K}$ is strictly convex.
Let $I$ be an interval of $\Z$ such that
$$
 I\cap [\bar i(y^{2m-1}),  \bar i(y^{2m-1})+\bar\nu(y^{2m-1})-1] =\varnothing
$$
for every odd iterate $y^{2m-1}$ with period larger than $aT$.
Then, for every $\epsilon>0$, we have
$$
 H_{S^1,k+e(K)} \bigl(\bar X^\epsilon(a,K),\bar X^{-\epsilon}(a,K);\Q\bigr) =0, \quad k\in\mathring I,
$$
where $\mathring I$ is $I$ minus its endpoints.
\end{theorem}

We will not prove Theorem \ref{thm:symmetric-Viterbo} in detail since when we replace $F_{a,K}$ in Section 7 of \cite{Vi89} by $\bar F_{a,K}$, the proof carries over verbatim after minor modification.
We only state the results needed in the proof of Theorem \ref{thm:main}.

By the homogeneity of $\bar F_{a,K}$, the inclusion $\bigl(\bar X^\epsilon(a,K),\bar X^{-\epsilon}(a,K)\bigr) \hookrightarrow \bigl(\bar X,\bar X^-(a,K)\bigr)$ is a weak equivariant homotopy equivalence (e.h.e for short).
Indeed, the functional $\bar F_{a,K}$ satisfies the Palais-Smale condition and has zero as its only critical level.
Hence $\bar X^\epsilon(a,K) \hookrightarrow \bar X$ is an e.h.e.
Since there is no critical value in $[-\epsilon,-\epsilon']$ for $0<\epsilon'<\epsilon$, $\bar X^{-\epsilon}(a,K) \hookrightarrow \bar X^{-\epsilon'}(a,K)$ is also an e.h.e.
Since $\bar X^-(a,K)=\bigcup_{\epsilon'>0}\bar X^{-\epsilon'}(a,K)$, the same conclusion holds for $\bar X^{-\epsilon}(a,K) \hookrightarrow \bar X^-(a,K)$.
Therefore, we obtain $H_{S^1,*} \bigl(\bar X^\epsilon(a,K),\bar X^{-\epsilon}(a,K);\Q\bigr) \cong H_{S^1,*} \bigl(\bar X,\bar X^-(a,K);\Q\bigr)$.

As in the 1st step of Theorem 7.1 of \cite{Vi89}, let $K'>K$ satisfy $\frac{K'T}{2\pi}\notin2\Z-1$.
Suppose $\bar H_{a,L}$ is strictly convex for every $L\in[K,K']$.
Then, for every $k\in\Z$, there is a natural isomorphism
\begin{equation}\label{eq:K-isomorphism}
 \xi_1: H_{S^1,k+e(K)}(\bar X,\bar X^-(a,K);\Q)
 \longrightarrow
 H_{S^1,k+e(K')}(\bar X,\bar X^-(a,K');\Q).
\end{equation}

We can rephrase the 2nd step of Theorem 7.1 \cite{Vi89} as follows.
For every finite set $\mathcal J\subset\Z$, there exists $a'>a$ with $a'T\notin\operatorname{per}_s(\Sigma)$ such that the map
\begin{equation}\label{eq:zero-composition}
 H_{S^1,k+e(K)}(\bar X,\bar X^-(a,K);\Q)
 \longrightarrow
 H_{S^1,k+e(K')}(\bar X,\bar X^-(a',K');\Q)
\end{equation}
is zero for every $k\in\mathcal J$ when $K'$ is sufficiently large and satisfies $K'T/(2\pi)\notin2\Z-1$.

Also, the 3rd step of Theorem 7.1 of \cite{Vi89} can be rephrased as follows.
Let $I$ be an interval of $\Z$.
Suppose that $ I\cap [\bar i(y^{2m-1}), \bar i(y^{2m-1})+\bar\nu(y^{2m-1})-1] =\varnothing $ for every odd iterate $y^{2m-1}$ with period in $[aT,a'T]$.
Then the map
$$
 H_{S^1,k+e(K')}(\bar X,\bar X^-(a,K');\Q)
 \longrightarrow
 H_{S^1,k+e(K')}(\bar X,\bar X^-(a',K');\Q)
$$
induced by the inclusion is an isomorphism for every $k\in\mathring I$.

By \eqref{eq:K-isomorphism} and \eqref{eq:zero-composition}, we obtain the following corollary, which will be used in the proof of Theorem \ref{thm:main}.
\begin{corollary}\label{cor:fixed-K-zero-map}
Let $a>0$ satisfy $aT\notin\operatorname{per}_s(\Sigma)$, and let $\mathcal J\subset\Z$ be finite.
There are $a'>a$ with $a'T\notin\operatorname{per}_s(\Sigma)$ and a constant $K_0>0$ such that, whenever $K'\geq K_0$ satisfies $\frac{K'T}{2\pi}\notin2\Z-1$ and $\bar H_{t,K'}$ is strictly convex for every $t\in[a,a']$, the homomorphism
$$
 j_k^{K'}(a,a'): H_{S^1,k+e(K')}(\bar X,\bar X^-(a,K');\Q)
 \longrightarrow
 H_{S^1,k+e(K')}(\bar X,\bar X^-(a',K');\Q)
$$
is zero for every $k\in\mathcal J$.
If $b>a'$ satisfies $bT\notin\operatorname{per}_s(\Sigma)$ and $\bar H_{t,K'}$ is strictly convex for every $t\in[a,b]$, then $j_k^{K'}(a,b)=0$ too.
\end{corollary}

We now prove the main theorem.
\begin{proof}[Proof of Theorem \ref{thm:main}]
We prove Theorem \ref{thm:main} by contradiction.
Assume $(\bf{F_s})$ and suppose that there exist some prime symmetric closed characteristics with zero symmetric mean indices.
Without loss of generality, we have the following condition.

$(\mathbf{C_s})$ Let $\{(\tau_k,y_k)\}_{1\leq k\leq q}$ be the prime symmetric closed characteristics on $\Sigma\in S\mathcal H_{\rm st}(4)$.
There exists an integer $q_0\in\{1,\ldots,q\}$ such that $\hat{\bar i}(y_k)=0$ for $1\leq k\leq q_0$ and $\hat{\bar i}(y_k)\neq0$ for $q_0+1\leq k\leq q$.

Assume $(\bf{C_s})$.
We use Lemma \ref{lem:R4-zero-mean-classification} and divide the proof into six steps.

\textbf{Step 1.}
Lemma \ref{lem:R4-zero-mean-classification} shows $\bar i(y_k^{2m-1})=-1$ for every $1\leq k\leq q_0$ and $m\in\N$.
Note that if $\hat{\bar i}(y_k)>0$ (respectively $<0$), then we must have $\bar i(y_k^{2m-1})\longrightarrow+\infty$ (respectively $-\infty$) as $m\to+\infty$.
Hence, for every $q_0+1\leq k\leq q$, we have $\bar i(y_k^{2m-1})\neq0,\ -1,\ -2$ whenever $m$ is sufficiently large.
Since there are only finitely many prime symmetric closed characteristics, we may choose $a$ large enough such that every odd iterate of $y_k$ with period larger than $aT$ has one of the following index properties:
\begin{equation}\label{eq:index-alternatives}
 \left\{
 \begin{array}{ll}
 \text{\rm either (i)}&
 \text{equal to }-1,\quad\text{when }\hat{\bar i}(y_k)=0 \quad \text{for } 1\le k\le q_0 \\[1mm]
 \text{\rm or (ii)}&
 \text{different from }0,-1\text{ and }-2,\quad
 \text{when }\hat{\bar i}(y_k)\neq0 \quad \text{for } q_0+1 \le k\le q .
 \end{array}
 \right.
\end{equation}

\textbf{Step 2.}
For $a\in \R$, let $\bar X^-(a,K)=\{x\in\bar X:\bar F_{a,K}(x)<0\}$.
From now on, we fix a large enough $a$ from \textbf{Step 1} so that $aT\notin\per_s(\Sigma)$.
Take $K=K(a)$ sufficiently large, with $\frac{KT}{2\pi}\notin2\Z-1$, so that $\bar H_{a,K}$ is strictly convex as in \eqref{eq:strict-convexity}.
Corollary \ref{cor:fixed-K-zero-map} shows that for $a'>a$ and sufficiently large $K'>K$ the inclusion maps in the degrees $0$ and $-1$ vanish.
We choose them so that $a'T\notin\per_s(\Sigma)$ and $\bar H_{t,K'}$ is strictly convex for all $t\in[a,a']$.
We may also require $(aT,a'T)$ to contain the period of at least one odd iterate of some $y_k$ with $1\leq k\leq q_0$.
Indeed, since these periods tend to $+\infty$, the last part of Corollary \ref{cor:fixed-K-zero-map} permits us to enlarge $a'$.

Now let $A=\bar X^-(a,K')$, $A'=\bar X^-(a',K')$.
For every $t\in[a,a']$, since $\bar F_{t,K'}$ is homogeneous of degree two, any critical orbit $S^1\cdot x$ of $\bar F_{t,K'}$ has critical value $\bar F_{t,K'}(S^1\cdot x)=0$.
Since $aT, a'T\notin\per_s(\Sigma)$, at $t=a$ and $t=a'$, the zero level contains no nonzero critical orbit, i.e., the boundary sets of $A$ and $A'$ contain no nonzero critical orbits.
Also, by definition the origin does not belong to $A$ or $A'$.
Hence
\begin{equation}\label{eq:sphere-isomorphism}
 H_{S^1,e(K')+i}(A',A;\Q) \cong
 H_{S^1,e(K')+i} (A'\cap S(\bar X),A\cap S(\bar X);\Q), \quad i\in\Z
\end{equation}
where $S(\bar X)$ is the unit sphere of $\bar X$.

Since each critical orbit $S^1\cdot x$ of $\bar F_{t,K'}|_{S(\bar X)}$ at level zero corresponds to an odd iterate $y_k^{2m-1}$ for some $m\in\N$ and $1\leq k\leq q$, we denote such an orbit by $S^1\cdot x_{t,k,m}$ and the period of the corresponding closed characteristic $y_k^{2m-1}$ by $T_{t,k,m}$.
From the construction of $\bar F_{t,K'}$, $T_{t,k,m}=tT$.
Thus every critical orbit in $(A'\setminus A)\cap S(\bar X)$ satisfies
$$
 aT\leq T_{t,k,m}\leq a'T.
$$
Since there are only finitely many critical orbits, there exist $\hat j\in\N$ and $a<t_1<\cdots<t_{\hat j}<a'$ such that $\bar F_{t,K'}$ has a critical orbit in $(A'\setminus A)\cap S(\bar X)$ only if $t\in\{t_1,\ldots,t_{\hat j}\}$.

Similar to \cite{DLLW24}, in order to compute the homology \eqref{eq:sphere-isomorphism}, we introduce below a new functional $\hat t$.
To define $\hat t$, we set
$$
 \mathcal U=\{x\in S(\bar X)\mid \bar F_{a,K'}(x)>0>\bar F_{a',K'}(x)\}.
$$

\begin{claim}\label{claim:derivative}
For every $(x,t)\in S(\bar X)\times[a,a']$, we have
$$
 \frac{\partial}{\partial t}\bar F_{t,K'}(x)<0.
$$
\end{claim}

In fact, by the definition of $\bar H_t(x)$ in Section~\ref{sec:4}, it is strictly increasing in $t$ when $x\neq0$, and then so is $\bar H_{t,K'}(x)$.
Then the Fenchel dual function $\bar H_{t,K'}^*(y)$ is strictly decreasing in $t$ when $y\neq0$.
Let
$$
 v_{t,K'}(x) =\nabla(\bar H_{t,K'})^*(-J\dot x+K'x).
$$
Then
$$
 \frac{\partial}{\partial t}\bar F_{t,K'}(x) =-\frac12\int_0^{T/2}j(v_{t,K'}(x)(s))^2\,ds.
$$
If the right-hand side is zero, then $v_{t,K'}(x)=0$ almost everywhere.
Hence $-J\dot x+K'x=0$. 
Since $\frac{K'T}{2\pi}\notin2\Z-1$, the operator $-J\frac d{ds}+K'$ is an isomorphism from $\bar X$ onto $\bar E$, and hence $x=0$. 
This is impossible on $S(\bar X)$.
Thus Claim \ref{claim:derivative} is proved.

For each $x\in\mathcal U$, Claim \ref{claim:derivative} implies that there is a unique $t\in(a,a')$ with $\bar F_{t,K'}(x)=0$.
By the implicit function theorem there is a $C^1$ map
$$
 \hat t:\mathcal U\longrightarrow(a,a')
$$
given by the equation
\begin{equation}\label{eq:barF}
 \bar F_{\hat t(x),K'}(x)=0,\quad x\in\mathcal U.
\end{equation}
Differentiating this identity in the $x$-variable gives
$$
 \frac{\partial}{\partial t}\bar F_{\hat t(x),K'}(x)
 \,\hat t'(x)+\bar F'_{\hat t(x),K'}(x)=0, \quad x\in\mathcal U.
$$
It follows that $x_0$ is a critical point of $\hat t$ with critical value $t_0=\hat t(x_0)$ if and only if it is a critical point of $\bar F_{t_0,K'}|_{S(\bar X)}$ with critical value 0.
Moreover, $\hat t$ is $S^1$-invariant.

\begin{claim}\label{claim:critical-module}
If $x_0$ is a critical point of $\hat t$ with critical value $t_0=\hat t(x_0)$, then the two critical modules are isomorphic:
\begin{equation}\label{eq:critical-module}
 C_{S^1,*}(\hat t,S^1\cdot x_0) \cong
 C_{S^1,*} (\bar F_{t_0,K'}|_{S(\bar X)},S^1\cdot x_0).
\end{equation}
\end{claim}

In fact, choose a small $S^1$-invariant neighborhood $U$ of $S^1\cdot x_0$ in $\mathcal U$.
By Claim \ref{claim:derivative} and \eqref{eq:barF}, $\hat t(x)\leq t_0$ if and only if $\bar F_{t_0,K'}(x)\leq0$ for all $x\in U$.
Therefore,
$$
 \left\{
 \begin{aligned}
 &\{x\in U\mid\hat t(x)\leq t_0\}
 =\{x\in U\mid \bar F_{t_0,K'}(x)\leq0\},\\
 &\{x\in U\mid\hat t(x)\leq t_0\}\setminus S^1\cdot x_0
 =\{x\in U\mid \bar F_{t_0,K'}(x)\leq0\}\setminus S^1\cdot x_0 .
 \end{aligned}
 \right.
$$
The definition of the $S^1$-critical module implies \eqref{eq:critical-module} in every degree. Thus Claim \ref{claim:critical-module} holds.

Choose $\epsilon>0$ so small that $a+\epsilon<a'-\epsilon$ and $cT\notin\per_s(\Sigma)$ for every $c\in[a,a+\epsilon]\cup[a'-\epsilon,a']$.
For $c\in[a,a']$, let $D_c=\{x\in S(\bar X)\mid\bar F_{c,K'}(x)\leq0\}$.
For $a<\kappa<a'$, let $\hat t^{\,\kappa}=\{x\in\mathcal U\mid\hat t(x)\leq\kappa\}$.
By Claim~\ref{claim:derivative} and \eqref{eq:barF}, we then obtain
$$
\left\{
\begin{aligned}
D_{a'-\epsilon}\cap\mathcal U
&=\{x\in\mathcal U\mid\bar F_{a'-\epsilon,K'}(x)\leq0\}=\hat t^{\,a'-\epsilon},\\
D_{a+\epsilon}\cap\mathcal U
&=\{x\in\mathcal U\mid\bar F_{a+\epsilon,K'}(x)\leq0\}=\hat t^{\,a+\epsilon}.
\end{aligned}
\right.
$$

For every $c\in[a,a+\epsilon]\cup[a'-\epsilon,a']$, $0$ is a regular value of $\bar F_{c,K'}|_{S(\bar X)}$.
Indeed, by homogeneity, a critical point of this restriction at level zero is a nonzero critical point of $\bar F_{c,K'}$, which implies $cT\in\per_s(\Sigma)$, contradicting the choice of $\epsilon$.
By equivariant deformation on $[a,a+\epsilon]$ and $[a'-\epsilon,a']$ and regularity of the zero levels (cf. \cite[Appendix~1, Proposition~3(a)]{Vi89}), the inclusions
$A\cap S(\bar X)\hookrightarrow D_{a+\epsilon}$ and $D_{a'-\epsilon}\hookrightarrow A'\cap S(\bar X)$ induce isomorphisms in equivariant homology.
Since $a+\epsilon<a'-\epsilon$, Claim~\ref{claim:derivative} gives $D_{a+\epsilon}\subset D_{a'-\epsilon}$.
We therefore have the following inclusions of pairs
$$
(A'\cap S(\bar X),A\cap S(\bar X))\hookleftarrow
(D_{a'-\epsilon},A\cap S(\bar X))\hookrightarrow
(D_{a'-\epsilon},D_{a+\epsilon}).
$$

Since $\bar F_{a',K'}<0$ on $D_{a'-\epsilon}$, we have
$$
D_{a'-\epsilon}\setminus\mathcal U=\{x\in D_{a'-\epsilon}\mid\bar F_{a,K'}(x)\leq0\}.
$$
The set $D_{a'-\epsilon}\setminus\mathcal U$ is closed in $D_{a'-\epsilon}$, and every $x\in D_{a'-\epsilon}\setminus\mathcal U$ has an open neighborhood in $D_{a'-\epsilon}$ contained in $D_{a+\epsilon}$ because $\bar F_{a+\epsilon,K'}$ is continuous and $\bar F_{a+\epsilon,K'}(x)<0$.
By equivariant excision, the inclusion
$$
(D_{a'-\epsilon}\cap\mathcal U,D_{a+\epsilon}\cap\mathcal U)
\hookrightarrow (D_{a'-\epsilon},D_{a+\epsilon})
$$
induces an isomorphism in relative equivariant homology.
Combining the above isomorphisms and using $D_c\cap\mathcal U=\hat t^{\,c}$, we obtain
\begin{align}
H_{S^1,e(K')+i}(A'\cap S(\bar X),A\cap S(\bar X);\Q)
&\cong H_{S^1,e(K')+i}(D_{a'-\epsilon},D_{a+\epsilon};\Q)\notag\\
&\cong H_{S^1,e(K')+i}(\hat t^{\,a'-\epsilon},\hat t^{\,a+\epsilon};\Q),\quad i\in\Z.
\label{eq:cong}
\end{align}

\textbf{Step 3.}
By \eqref{eq:index-alternatives}, no odd iterate of $y_k$ with period in $[aT,a'T]$ has symmetric index $0$ or $-2$.
Therefore, by the discussion on pp. 78-79 of \cite{Ch93}, we obtain
\begin{align*}
 C_{S^1,e(K')} \bigl(\bar F_{t_j,K'}|_{S(\bar X)},S^1\cdot x_{t_j,k,m}\bigr)
 =C_{S^1,e(K')-2} \bigl(\bar F_{t_j,K'}|_{S(\bar X)},S^1\cdot x_{t_j,k,m}\bigr)=0.
\end{align*}
By \eqref{eq:critical-module}, the corresponding critical modules of $\hat t$ also vanish, we have
\begin{equation}\label{eq:local-adjacent-zero}
 \begin{aligned}
 C_{S^1,e(K')}(\hat t,S^1\cdot x_{t_j,k,m})
 = C_{S^1,e(K')-2}(\hat t,S^1\cdot x_{t_j,k,m}) =0.
 \end{aligned}
\end{equation}
The gradient formula for $\bar F_{t,K'}$, together with homogeneity and the compact embedding $\bar X\hookrightarrow\bar E$, gives the Palais-Smale condition for $\hat t$ on $\hat t^{-1}([a+\epsilon,a'-\epsilon])$.
Combining \eqref{eq:local-adjacent-zero} with equivariant Morse inequalities, we have
\begin{equation}\label{eq:H-adjacent-zero}
 H_{S^1,e(K')} (\hat t^{\,a'-\epsilon},\hat t^{\,a+\epsilon};\Q)
 = H_{S^1,e(K')-2} (\hat t^{\,a'-\epsilon},\hat t^{\,a+\epsilon};\Q) =0.
\end{equation}
Therefore, combining \eqref{eq:H-adjacent-zero} with \eqref{eq:sphere-isomorphism} and \eqref{eq:cong}, we obtain
\begin{equation}\label{eq:HAA-zero}
 H_{S^1,e(K')}(A',A;\Q) = H_{S^1,e(K')-2}(A',A;\Q) =0.
\end{equation}

\textbf{Step 4.}
Consider the long exact sequence of the triple $(\bar X,A',A)$:
\begin{equation}\label{eq:long-exact}
\begin{aligned}
 \longrightarrow\;&
 H_{S^1,e(K')}(A',A) \xrightarrow{i_{0*}}
 H_{S^1,e(K')}(\bar X,A) \xrightarrow{j_{0*}}
 H_{S^1,e(K')}(\bar X,A')\\
 \xrightarrow{\partial_{0*}}\;&
 H_{S^1,e(K')-1}(A',A) \xrightarrow{i_{-1*}}
 H_{S^1,e(K')-1}(\bar X,A) \xrightarrow{j_{-1*}}
 H_{S^1,e(K')-1}(\bar X,A')\\
 \xrightarrow{\partial_{-1*}}\;&
 H_{S^1,e(K')-2}(A',A) \longrightarrow \cdots.
\end{aligned}
\end{equation}
All homology groups in \eqref{eq:long-exact} have coefficients in $\Q$.
Now consider the composition
$$
 H_{S^1,e(K)}(\bar X,\bar X^-(a,K);\Q)
 \xrightarrow{\ \xi_1\ }
 H_{S^1,e(K')}(\bar X,\bar X^-(a,K');\Q)
 \xrightarrow{\ \xi_2\ }
 H_{S^1,e(K')}(\bar X,\bar X^-(a',K');\Q),
$$
where $\xi_1$ is the isomorphism by \eqref{eq:K-isomorphism} and $\xi_2=j_{0*}$ by \eqref{eq:long-exact}.
Their composition $\xi_2 \circ \xi_1$ is the zero map by \eqref{eq:zero-composition}.
Therefore $j_{0*}=0$.

Combining \eqref{eq:HAA-zero} and \eqref{eq:long-exact} we obtain
\begin{align}
 H_{S^1,e(K')}(\bar X,A;\Q) =\ker j_{0*} =\operatorname{im}i_{0*}
 =i_{0*}\bigl(H_{S^1,e(K')}(A',A;\Q)\bigr) =0.                                         
\label{eq:absolute-A-zero}
\end{align}

We now fix $a'$ and apply Corollary~\ref{cor:fixed-K-zero-map} with $a'$ and $\mathcal J=\{0,-1\}$ to choose $a''>a'$ with $a''T\notin\operatorname{per}_s(\Sigma)$.
We enlarge $K'$ if necessary so that all the arguments from \textbf{Step 1} to \textbf{Step 4} are valid when we replace $(a,a')$ by $(a',a'')$.
Let $A''=\bar X^-(a'',K')$. Applying the same argument to the triple $(\bar X,A'',A')$, we obtain
\begin{equation}\label{eq:XA-zero}
 H_{S^1,e(K')}(\bar X,A';\Q)=0.
\end{equation}
Combining \eqref{eq:HAA-zero}, \eqref{eq:absolute-A-zero} and \eqref{eq:XA-zero} with the long exact sequence \eqref{eq:long-exact}, we obtain the short exact sequence
\begin{equation}\label{eq:short-exact}
 0\longrightarrow
 H_{S^1,e(K')-1}(A',A;\Q) \xrightarrow{i_{-1*}}
 H_{S^1,e(K')-1}(\bar X,A;\Q) \xrightarrow{j_{-1*}}
 H_{S^1,e(K')-1}(\bar X,A';\Q) \longrightarrow0.
\end{equation}

\textbf{Step 5.}
For $1\leq k\leq q_0$, Lemma \ref{lem:R4-zero-mean-classification} implies that there are infinitely many odd iterates of $y_k$ of symmetric index $-1$.
For the above chosen large enough $a'>a$, only finitely many of these iterates have periods in the interval $[aT,a'T]$.
Therefore, their corresponding critical orbits $S^1\cdot x_{t_j,k,m}$ all have Morse index $e(K')-1$.
By the equivariant version of Theorem~I.4.2 and the discussion on pp. 78-79 of \cite{Ch93}, we obtain
\begin{align}
 H_{S^1,e(K')-1}(A',A;\Q)
 &\cong H_{S^1,e(K')-1} (A'\cap S(\bar X),A\cap S(\bar X);\Q)\notag\\
 &\cong H_{S^1,e(K')-1} (\hat t^{\,a'-\epsilon},\hat t^{\,a+\epsilon};\Q)\notag\\
 &\cong\bigoplus_{\substack{
       aT\leq T_{t_j,k,m}\leq a'T\\
       1\leq k\leq q_0}}
 C_{S^1,e(K')-1} (\hat t,S^1\cdot x_{t_j,k,m})\notag\\
 &\cong\bigoplus_{\substack{
       aT\leq T_{t_j,k,m}\leq a'T\\
       1\leq k\leq q_0}}
 C_{S^1,e(K')-1} (\bar F_{t_j,K'}|_{S(\bar X)},S^1\cdot x_{t_j,k,m})\notag\\
 &\cong\bigoplus_{\substack{
       aT\leq T_{t_j,k,m}\leq a'T\\
       1\leq k\leq q_0}}\Q
 \neq0.           
\label{eq:Oplus}
\end{align}
The first two isomorphisms follow from \eqref{eq:sphere-isomorphism} and \eqref{eq:cong} respectively.
The fourth one follows from \eqref{eq:critical-module}.
We explain the third isomorphism as follows.
Denote by
$$
\begin{aligned}
 M_l(a,a')
 &=
 \sum_{\substack{
       aT\leq T_{t_j,k,m}\leq a'T\\
       1\leq k\leq q}}
 \operatorname{rank}
 C_{S^1,l}(\hat t,S^1\cdot x_{t_j,k,m}),\\
 \beta_l(a,a')
 &= \operatorname{rank}
 H_{S^1,l}(\hat t^{\,a'-\epsilon},  \hat t^{\,a+\epsilon};\Q).
\end{aligned}
$$
By \eqref{eq:local-adjacent-zero},
$$
 M_{e(K')}(a,a')=M_{e(K')-2}(a,a')=0,
$$
whereas \eqref{eq:H-adjacent-zero} gives
$$
 \beta_{e(K')}(a,a')=\beta_{e(K')-2}(a,a')=0.
$$
These two equations and the equivariant Morse inequalities imply
$$
 M_{e(K')-1}(a,a')=\beta_{e(K')-1}(a,a').
$$
Then the third isomorphism holds.

\textbf{Step 6.}
By the exactness of the sequence \eqref{eq:short-exact} and \eqref{eq:Oplus}, we obtain
$$
\begin{aligned}
 H_{S^1,e(K')-1}(\bar X,A;\Q)
 &\cong H_{S^1,e(K')-1}(A',A;\Q) \oplus H_{S^1,e(K')-1}(\bar X,A';\Q) \neq0.
\end{aligned}
$$

We may also require $(a'T,a''T)$ to contain the period of at least one odd iterate of some $y_k$ with $1\leq k\leq q_0$.
By Corollary~\ref{cor:fixed-K-zero-map}, 
choose $b>a''$ with $bT\notin\operatorname{per}_s(\Sigma)$ and enlarge $K'$ if necessary so that the inclusion $A''\subset\bar X^-(b,K')$ induces the zero map $H_{S^1,e(K')}(\bar X,A'';\Q)\rightarrow H_{S^1,e(K')}(\bar X,\bar X^-(b,K');\Q)$.
We use the same $K'$ in the preceding arguments and in the definitions of $A$, $A'$, and $A''$, 
choosing it sufficiently large for the applications of Corollary~\ref{cor:fixed-K-zero-map} to $[a,a']$, $[a',a'']$, 
and $[a'',b]$, with $K'T/(2\pi)\notin2\Z-1$ and $\bar H_{t,K'}$ strictly convex for every $t\in[a,b]$.
The argument of \textbf{Step~4} applied to $(a'',b)$ gives $H_{S^1,e(K')}(\bar X,A'';\Q)=0$.
By \textbf{Step~5} applied to $(a',a'')$, we have $H_{S^1,e(K')-1}(A'',A';\Q)\neq0$.
The long exact sequence of the triple $(\bar X,A'',A')$ therefore yields
\begin{equation}\label{eq:HXA-nonzero}
 H_{S^1,e(K')-1}(\bar X,A';\Q)\neq0.
\end{equation}

Finally, we can get a contradiction by computing $j_{-1*}$ in two different ways.
If $j_{-1*}=0$, the exactness of \eqref{eq:short-exact} implies
$$
 H_{S^1,e(K')-1}(\bar X,A';\Q) =\ker\partial_{-1*} =\operatorname{im}j_{-1*} =0,
$$
which contradicts \eqref{eq:HXA-nonzero}.
Therefore $j_{-1*}\neq0$.
However, by the same arguments of the fact that $j_{0*}=0$ in \textbf{Step 4}, we know that $j_{-1*}=0$, which yields a contradiction.
The proof of Theorem \ref{thm:main} is completed.

\end{proof}

We now return to the variational structure of Section \ref{sec:2} to prove the stability result.

\begin{proof}[Proof of Theorem \ref{thm:two-symmetric-elliptic}]
Let $(\tau_1,y_1)$ and $(\tau_2,y_2)$ be the two geometrically distinct prime symmetric closed characteristics on $\Sigma$.
For $j=1,2$, let $\gamma_j:[0,\tau_j/2]\to\Sp(4)$ be the symplectic path associated with $y_j$.
By Lemma \ref{lem:first-block}, $\gamma_j(\tau_j/2)$ is conjugate to $N_1(-1,-1)\diamond M_j$ for some $M_j\in\Sp(2)$.
By symmetry, the symplectic path associated with $y_j$ on $[0,\tau_j]$ is $\gamma_j^2$.
Thus $\gamma_j^2(\tau_j)=[\gamma_j(\tau_j/2)]^2$ is conjugate to $N_1(-1,-1)^2\diamond M_j^2$.
Since every iterate of $y_j$ is non-degenerate, $1$ is not an eigenvalue of $M_j^{2m}$ for any $m\in\N$.
In particular, no eigenvalue of $M_j$ is a root of unity.
Therefore $\ker(M_j^{4m-2}-I_2)=\{0\}$ for every $m\in\N$.
Since $M_j^{4m-2}-I_2=(M_j^{2m-1}-I_2)(M_j^{2m-1}+I_2)$, we have 
$$
\ker(M_j^{2m-1}+I_2) \subseteq\ker(M_j^{4m-2}-I_2)=\{0\}.
$$
Therefore by \eqref{eq:symmetric-index-nullity} and \eqref{eq:half-period-decomposition}, for every $m\in\N$, we have 
\begin{align*}
 \bar\nu(y_j^{2m-1})
 &=\dim\ker\bigl([\gamma_j(\tau_j/2)]^{2m-1}+I_4\bigr)\\
 &=\dim\ker\bigl(N_1(-1,-1)^{2m-1}+I_2\bigr)
   +\dim\ker\bigl(M_j^{2m-1}+I_2\bigr)\\
 &=1.
\end{align*}
Thus Theorem \ref{thm:main} gives $\hat{\bar i}(y_j)\neq0$ for $j=1,2$.

By \eqref{eq:symmetric-average-euler}, the average Euler characteristic $\hat{\bar\chi}(y_j)$ of $y_j$ satisfies $\hat{\bar\chi}(y_j)=(-1)^{\bar i(y_j)}$.
Then the resonance identities \eqref{eq:symmetric-resonance} in Theorem~\ref{thm:symmetric-resonance} give
\begin{equation}\label{eq:two-symmetric-resonance}
 \sum_{\substack{1\leq j\leq 2\\ \hat{\bar i}(y_j)>0}}
 \frac{(-1)^{\bar i(y_j)}}{\hat{\bar i}(y_j)}=1, \quad 
 \sum_{\substack{1\leq j\leq 2\\ \hat{\bar i}(y_j)<0}}
 \frac{(-1)^{\bar i(y_j)}}{\hat{\bar i}(y_j)}=0.
\end{equation}
By the first identity in \eqref{eq:two-symmetric-resonance}, at least one of the two symmetric mean indices is positive.
If the other one is negative, then the sum in the second identity contains one nonzero term and is therefore nonzero, which contradicts \eqref{eq:two-symmetric-resonance}.
Hence $\hat{\bar i}(y_j)>0$ for $j=1,2$, and
\begin{equation}\label{eq:two-symmetric-positive-resonance}
 \frac{(-1)^{\bar i(y_1)}}{\hat{\bar i}(y_1)} +\frac{(-1)^{\bar i(y_2)}}{\hat{\bar i}(y_2)}=1.
\end{equation}

Since $N_1(-1,-1)^2=N_1(1,2)$ is symplectically conjugate to $N_1(1,1)$ and $M_j\in\Sp(2)$ has no eigenvalue that is a root of unity, the above decomposition of $\gamma_j^2(\tau_j)$ shows that $y_j$ is either hyperbolic or irrationally elliptic.
If $y_j$ is hyperbolic, then $p_-=p_0=p_+=r=0$ and $s=1$ in the normal form of $M_j$ in \eqref{eq:M-normal-form}, so \eqref{eq:R4-precise-odd-formula} and \eqref{eq:symmetric-mean-formula} give
\begin{equation}\label{eq:stability-hyperbolic-iteration}
 \hat{\bar i}(y_j)=i_1(\gamma_j)\in\Z,  \quad
 \bar i(y_j^{2m-1})=(2m-1)\hat{\bar i}(y_j)-1, \quad m\in\N.
\end{equation}
If $y_j$ is irrationally elliptic, then $M_j$ is conjugate to $R(\theta_j)$ with $\theta_j/\pi\notin\Q$.
In this case, \eqref{eq:symmetric-mean-formula} gives
$$
 \hat{\bar i}(y_j)=i_1(\gamma_j)-1+\frac{\theta_j}{\pi}\notin\Q.
$$

Suppose that one of the two symmetric closed characteristics is hyperbolic, and denote it by $y_1$.
Since $\hat{\bar i}(y_1) \in \Z$, \eqref{eq:two-symmetric-positive-resonance} implies that $\hat{\bar i}(y_2) \in \Q$.
Thus $y_2$ is also hyperbolic.

By \eqref{eq:stability-hyperbolic-iteration} and $\hat{\bar i}(y_j)>0$ for $j=1,2$, we have $\hat{\bar i}(y_j)\in\N$ and $\bar i(y_j^{2m-1})=(2m-1)\hat{\bar i}(y_j)-1\geq2m-2\geq0$ for $j=1,2$ and $m\in\N$.
The symmetric mean indices are nonzero, so \eqref{eq:symmetric-morse-series} applies.
For each fixed $p\in\Z$ and sufficiently large $a$, the Morse type number $M_p$ in \eqref{eq:symmetric-morse-series} is given by \eqref{eq:normalized-symmetric-morse-series} as
$$
 M_p= \sum_{\substack{1\leq j\leq 2,\ m\in\N\\ (2m-1)\tau_j<aT}}
 \dim C_{S^1,e(K)+p} \bigl(\bar F_{a,K},S^1\cdot x_j^{2m-1}\bigr),
$$
where $x_j^{2m-1}$ is the critical point of $\bar F_{a,K}$ corresponding to $y_j^{2m-1}$.
By the isomorphism between the corresponding critical modules of $\bar F_{a,K}$ and $\bar\Psi_{a,K}$, it follows from \eqref{eq:local-module} that $y_j^{2m-1}$ contributes $1$ to $M_p$ if $\bar i(y_j^{2m-1})=p$ and contributes $0$ to $M_p$ if $\bar i(y_j^{2m-1})\neq p$.
The Betti numbers from Section~\ref{sec:2} satisfy
$$
 \beta_p=\dim H_p(\CP^\infty;\Q)
  =\begin{cases}
  1,&p\in2\N^0,\\
  0,&p\notin2\N^0.
 \end{cases}
$$
By \eqref{eq:symmetric-morse-coefficients}, we have $M_0\geq\beta_0=1$, so $\bar i(y_{j_0}^{2m_0-1})=0$ for some $j_0\in\{1,2\}$ and $m_0\in\N$.
By \eqref{eq:stability-hyperbolic-iteration}, we have $(2m_0-1)\hat{\bar i}(y_{j_0})=1$.
Since $\hat{\bar i}(y_{j_0})\in\N$, we have $m_0=1$ and $\hat{\bar i}(y_{j_0})=1$.
Without loss of generality, we may assume that $\hat{\bar i}(y_1)=1$.
Then by \eqref{eq:stability-hyperbolic-iteration}, we have
$$
 \bar i(y_1^{2m-1})=2m-2,\quad m\in\N.
$$
Thus the odd iterates of $y_1$ contribute $1$ to $M_p$ for each $p\in2\N^0$ and contribute $0$ to $M_p$ for each $p\notin 2\N^0$, so their contribution is exactly $\beta_p$.

By \eqref{eq:stability-hyperbolic-iteration}, we also have
$$
 \bar i(y_2^3)=3\bar i(y_2)+2>\bar i(y_2)+1.
$$
Since $\bar i(y_2^{2m-1})$ increases with $m$, the only odd iterate of $y_2$ whose symmetric index lies in $[0,\bar i(y_2)+1]$ is $y_2$ itself.
Thus, for $0\leq p\leq\bar i(y_2)+1$, we have $M_p-\beta_p=1$ if $p=\bar i(y_2)$ and $M_p-\beta_p=0$ otherwise.
Since all odd iterates have nonnegative symmetric indices, we have $M_p=0$ for $p<0$.
Taking $k=\bar i(y_2)+1$ in \eqref{eq:stability-Morse}, we obtain
$$
 0\leq\sum_{p=0}^{\bar i(y_2)+1}(-1)^{\bar i(y_2)+1-p}(M_p-\beta_p)=-1,
$$
which is a contradiction.
Therefore both $y_1$ and $y_2$ are irrationally elliptic.
\end{proof}

\bibliographystyle{abbrv}

\end{document}